\documentclass[a4paper,11pt]{amsart}

\usepackage{amssymb,amsmath,eufrak}
\author[Florent Benaych-Georges]{Florent Benaych-Georges}\address{Florent Benaych-Georges, LPMA,  UPMC Univ Paris 6, Case courier 188, 4, Place Jussieu, 75252 Paris Cedex 05, France} \email{florent.benaych@gmail.com}
\title{On a surprising relation between rectangular and square free convolutions}
\date{\today}
\newcommand{\inv}{^{\lan -1\ran}}
\newcommand{\be}{\begin{equation}}
\newcommand{\ee}{\end{equation}}
\newcommand{\beq}{\begin{eqnarray*}}
\newcommand{\eeq}{\end{eqnarray*}}

\newcommand{\lan}{\langle}
\newcommand{\ran}{\rangle}

\newcommand{\R}{\mathbb{R}}
\newcommand{\C}{\mathbb{C}}

\newcommand{\ud}{\mathrm{d}}

\newcommand{\pro}{probability }

\newcommand{\f}{\frac}
\newcommand{\ff}{\frac{1}}
\newcommand{\lf}{\left}
\newcommand{\ri}{\right}

\newcommand{\st}{such that }
\newcommand{\la}{\lambda}

\newcommand{\vfi}{\varphi}

\newcommand{\mc}{\mathcal }

\newcommand{\eps}{\varepsilon}
\newcommand{\arc}{\boxplus_{\la}}

\newcommand{\bxp}{\boxplus}
\newcommand{\bxt}{\boxtimes}

\newcommand{\A}{\mc{A}}

\newcommand{\bxs}{\begin{picture}(9,5)  \put(2,0){\framebox(5,5){$\smallsetminus$}} \end{picture}}

\newcommand{\bck}{\backslash}
\newcommand{\esm}{\operatorname{ESM}}

\newtheorem{Th}{Theorem}
\newtheorem{propo}[Th]{Proposition}

\newtheorem{rmq}[Th]{Remark}
\newtheorem{cor}[Th]{Corollary}

\newenvironment{pr}{\noindent {\bf Proof. }}{\ \ \ $\square$}

\long\def\symbolfootnote[#1]#2{\begingroup
\def\thefootnote{\fnsymbol{footnote}}\footnote[#1]{#2}\endgroup}

\begin{document}
\maketitle

\symbolfootnote[0]{{\bf MSC 2000 subject classifications.}  
46L54, %Free proabability and free factors
15A52} %Random matrices 

\symbolfootnote[0]{{\bf Key words.} free probability,  random matrices, free convolution} 

\begin{abstract}Debbah and Ryan have recently \cite{dr07} proved a result about the limit empirical singular distribution of the sum of two rectangular random matrices whose dimensions tend to infinity. In this paper, we reformulate it in terms of the rectangular free convolution introduced in \cite{bg07} and then we give a new, shorter, proof of this result under weaker hypothesis: we do not suppose the \pro measure in question in this result to be compactly supported anymore. At last, we discuss the inclusion of this result in the family of relations between rectangular and square random matrices.\end{abstract}
%arXiv abtract: Debbah and Ryan have recently   proved a result about the limit empirical singular distribution of the sum of two rectangular random matrices as the dimensions tend to infinity. In this paper, we reformulate it in terms of the rectangular free convolution introduced in a previous paper  and then we give a new, shorter, proof of this result under weaker hypothesis: we do not suppose the \pro measure in question in this result to be compactly supported anymore. At last, we discuss the inclusion of this result in the family of relations between rectangular and square random matrices.

%\tableofcontents

\section*{Introduction}Free convolutions are operations on \pro measures on the real line which allow to compute the spectral or singular empirical measures of large random matrices which are expressed as sums or products of independent random matrices, the spectral measures of which are known. 
More specifically, the operations $\bxp,\bxt$, called respectively {\em free additive and multiplicative convolutions} are defined in the following way \cite{vdn91}. Let, for each $n$, $M_n$, $N_n$ be $n$ by $n$ independent random hermitian matrices, one of them having a distribution which is invariant under the action of the unitary group by conjugation, which empirical spectral measures\footnote{The {\em empirical spectral measure} of a matrix is the uniform law on its eigenvalues with multiplicity.}  converge, as $n$ tends to infinity, to non random \pro measures denoted respectively  by $\tau_1, \tau_2$. Then $\tau_1\bxp\tau_2$ is the limit of the empirical spectral law of $M_n+N_n$ and, in the case where the matrices are positive, $\tau_1\bxt\tau_2$ is the limit of the empirical spectral law of $M_nN_n$. In the same way, for any $\la\in [0,1]$, the {\em rectangular free convolution} $\arc$ is defined, in  \cite{bg07}, in the following way. Let $M_{n,p}, N_{n,p}$ be $n$ by $p$ independent random  matrices, one of them having a distribution which is   invariant by multiplication by any unitary matrix on any side,    which symmetrized\footnote{The {\em symmetrization} of a \pro measure $\mu$ on $[0,+\infty)$ is the law of $\eps X$, for $\eps, X$ independent random variables with respective laws $\f{\delta_1+\delta_{-1}}{2}, \mu$. Dealing with laws on $[0,+\infty)$ or with their symmetrizations is equivalent, but for historical reasons, the rectangular free convolutions have been defined with symmetric laws. In all this paper, we shall often pass from symmetric \pro measures to measures on $[0,+\infty)$ and vice-versa. Thus in order to avoid confusion, we shall mainly use the letter $\mu$ for measures on $[0,\infty)$ and $\nu$ for symmetric ones.} empirical singular measures\footnote{The {\em empirical singular measure} of a matrix $M$ with size $n$ by $p$ ($n\leq p$) is the 
empirical spectral measure of $|M|:=\sqrt{MM^*}$.} tend, as $n,p$ tend to infinity in such a way that $n/p$ tends to $\la$,  to non random \pro measures $\nu_1,\nu_2$. % which symmetrized square roots\footnote{The {\em symmetrized square root} of a law $\rho$ on $[0,+\infty)$ is the law of $\eps \sqrt{X}$, for $\eps, X$ independent random variables with respective laws $\f{\delta_1+\delta_{-1}}{2}, \rho$. We shall denote it by $\sqrt{\rho}$: it is the unique symmetric law which push-forward by the map $t\mapsto t^2$ is $\rho$. Considering laws on $[0,+\infty)$ or their symmetrized square roots is equivalent, but for historical reasons, the rectangular free convolutions have been defined with symmetric laws.}, are $\mu,\nu$. 
Then the symmetrized  empirical singular law of $M_{n,p}+N_{n,p}$ tends to  $\nu_1\arc \nu_2$.

These operations can be explicitly computed using either a combinatorial or an analytic machinery (see \cite{vdn91} and \cite{ns06} for $\bxp, \bxt$ and  \cite{bg07} for $\arc$). In the cases $\la=0$ or  $\la=1$, i.e. where the rectangular random matrices we consider are either ``almost flat" or ``almost square",  the rectangular free convolution with ratio $\la$ can be expressed with the additive free convolution:  $\bxp_1=\bxp$ and  for all symmetric laws $\nu_1,\nu_2$, $\nu_1\bxp_0 \nu_2$ is the symmetric law which push-forward by    the map $t\mapsto t^2$ is the free convolution of the push forwards of $\nu_1$ and $\nu_2$ by the same map. However, though one can find many analogies between  the definitions of $\bxp$ and $\arc$ and still more analogies have been proved \cite{fbg05.inf.div}, 
no general relation 
between $\arc$ and $\bxp$ had been proved until a paper of Debbah and Ryan \cite{dr07} (which submitted version, more focused on applications than on this result, is \cite{dr08}). It is to notice  that this result is not  due to researchers from the communities of Operator Algebras or Probability Theory, 
but to researchers from Information Theory, working on communication networks. 
In \cite{dr07}, Debbah and Ryan proved a result about random matrices which can be  interpreted as an expression, for certain \pro measures $\nu_1,\nu_2$, of their rectangular convolution  $\nu_1\arc\nu_2$ in terms of $\bxp$ and of another convolution, called the {\em free multiplicative deconvolution} and denoted by $\bxs$. In this  note, we present  this result with a new approach   and we give a new and  shorter proof, where the hypothesis are more general. This generalization of the hypothesis answers a question asked by Debbah and Ryan in the last section of their  paper \cite{dr07}. The question of a more general relation between square and rectangular free convolutions is considered in a last ``perspectives" section.

{\bf Acknowledgments:} The author would like to thank Raj Rao for bringing the paper \cite{dr07} to his attention and M\'erouane Debbah for his encouragements and many useful discussions. 

\section{The result of Debbah and Ryan}Let us define the operation $\bxs$ on certain pairs of \pro measures on $[0,+\infty)$ in the following way. For $\mu,\mu_2$ \pro measures on $[0,+\infty)$, if there is a \pro measure on $[0,+\infty)$ \st $\mu=\mu_1\bxt\mu_2$, then $\mu_1$ is called  the  {\em free multiplicative deconvolution} of $\mu$ by $\mu_2$ and is denoted by $\mu_1=\mu\bxs\mu_2$.
Let us define, for $\la\in (0,1]$, $\mu_\la$ be the law of $\la X$ for $X$ random variable distributed according to the Marchenko-Pastur law with parameter $1/\la$, i.e. the law with support $[(1-\sqrt{\la})^2, (1+\sqrt{\la})^2]$ and density $$x\mapsto \f{\sqrt{4\la -(x-1-\la)^2}}{2\pi\la x}.$$

Theorem 1 of \cite{dr07} states the following result. $\la\in (0,1]$ is fixed and $(p_n)$ is a sequence of positive integers \st $n/p_n$ tends to $\la$ as $n$ tends to infinity. $\delta_1$ denotes the Dirac mass at $1$.
\begin{Th}[Debbah and Ryan]\label{1.07.08.3}Let, for each $n$, $A_n$, $G_n$ be independent $n$ by $p_n$ random matrices \st the empirical spectral law of $A_nA_n^*$ converges almost surely weakly, as $n$ tends to infinity, to a compactly supported  \pro measure $\mu_A$ and \st the entries of $G_n$ are independent $N(0, \ff{p_n})$ random variables. Then the empirical spectral law of $(A_n+G_n)(A_n+G_n)^*$ converges almost surely to a compactly supported \pro measure   $\rho$ which, in the case where $\mu_A\bxs\mu_\la$ exists, satisfies the relation \be\label{1.07.08.2}\rho=[(\mu_A\bxs\mu_\la)\bxp\delta_1]\bxt\mu_\la.\ee
\end{Th} 

\begin{rmq}\label{1.7.8.16h}{\rm Note that in the case where $\mu_A\bxs\mu_\la$ doesn't exist, the relation \eqref{1.07.08.2} stays true in the formal sense. More specifically, for $\mu_A$ \pro measure \st $\mu_A\bxs\mu_\la$ exists, the moments of $\rho=[(\mu_A\bxs\mu_\la)\bxp\delta_1]\bxt\mu_\la$ have a polynomial expression in  the moments of $\mu_A$ (this can easily be seen by the theory of free cumulants \cite{ns06}). It happens that this relation between the moments of the limit  spectral law $  \rho$ of $(A_n+G_n)(A_n+G_n)^*$  and the ones of $\mu_A$ stays true even when $\mu_A\bxs\mu_\la$ doesn't exist. It follows from the original proof of Theorem \ref{1.07.08.3} and it will also follow from our proof (see Remark \ref{1.7.8.1}).}\end{rmq}

Note that by the very definition of the rectangular free convolution $\arc$ with ratio $\la$ recalled in the introduction and since the limit empirical spectral law of $GG^*$ is $\mu_\la$ (it is a well known fact, see, e.g. Theorem 4.1.9 of \cite{hiai}), this result can be stated as  follows: for all compactly supported \pro measure $\mu$ on $[0,+\infty)$ \st $\mu\bxs\mu_\la$ exists, \be\label{30.06.08.1}(\sqrt{\mu}\arc \sqrt{\mu_\la})^2=[ (\mu\bxs\mu_\la)\bxp\delta_1]\bxt\mu_\la,\ee where for any \pro measure $\rho$ on $[0,+\infty)$, $\sqrt{\rho}$ denotes the symmetrization of the push-forward by the  square root functions  of $\rho$ and for any symmetric \pro measure $\nu$ on the real line, $\nu^2$ denotes the push-forward of $\nu$ by the function $t\mapsto t^2$. This formula allows to express the operator $\arc\sqrt{\mu_\la}$ on the set of symmetric compactly supported \pro measures on the real line in terms of $\bxp$ and $\bxt$: for all symmetric \pro measure on the real line $\nu$, \be\label{30.06.08.2}\nu\arc \sqrt{\mu_\la}=\sqrt{ [(\nu^2\bxs\mu_\la)\bxp\delta_1]\bxt\mu_\la}.\ee

\section{A proof  of the generalized theorem of Debbah and Ryan}
$\la\in (0,1]$ is still fixed.
In this section, we shall give a new shorter proof  of the theorem of Debbah and Ryan, under weaker hypothesis: we shall prove  \eqref{30.06.08.2} without supposing the support of $\nu$ to be compactly supported. The proof   is based on the machinery of the rectangular free convolution and of the rectangular $R$-transform.

\subsection{Some analytic transforms} 
Let us first recall a few  facts about the analytic approach to $\bxt$ and $\arc$. Let us define, for $\rho$
\pro measure on $[0,\infty)$, $$M_\rho(z):=\int_{t\in \R}\f{zt}{1-zt}\ud \rho(t),\quad S_\rho(s)=\f{1+z}{z}M_\rho^{\lan -1\ran}(z),$$where, as it shall be in the rest of the text, the exponent $\inv$ stands for the inversion of analytic  functions on $\C\bck [0,+\infty)$ with respect to the composition operation $\circ$, in a neighborhood of zero. By \cite{vdn91}, for all pair $\mu_1, \mu_2$ of \pro measures on $[0,+\infty)$, $\mu_1\bxt\mu_2$ is characterized by the fact that $S_{\mu_1\bxt\mu_2}=S_{\mu_1}S_{\mu_2}$.

In the same way, the rectangular free convolution with ratio $\la$ can be computed with an analytic transform of \pro measures. Let $\nu$ be a symmetric \pro measure on the real line. Let us define 
$H_\nu(z)= z(\la M_{\nu^2}(z)+1)(M_{\nu^2}(z)+1)$. Then the {\em rectangular $R$-transform with ratio $\la$} of $\nu$ is defined to be $$C_\nu(z)=U\lf( \f{z}{H_\nu^{\lan -1\ran}(z)}-1\ri), $$where $U(z)=  \f{-\la-1+\lf[(\la+1)^2+4\la z\ri]^{1/2}}{2\la}$. By theorem 3.12 of  \cite{bg07}, for all pair $\nu_1, \nu_2$ of symmetric \pro measures, $\nu_1\arc\nu_2$ is characterized by the fact that $C_{\nu_1\arc\nu_2}=C_{\nu_1}+C_{\nu_2}$.  

\subsection{Some preliminary computations}
Note that by \cite{ns06}, the $S$- and $R$-transforms of a \pro measure $\mu$ on $[0,+\infty)$  are linked by the relation $S_\mu(s)=\ff{z}R_\mu\inv(z)$, thus since the free cumulants of the Marchenko-Pastur law with parameter $1/\la$ are all equal to $1/\la$ (see \cite{ns06}), we have  $S_{\mu_\la}(z)=\ff{1+\la z}$. Moreover, since by \cite{ns06} again, $S_\mu(s)= \f{1+z}{z}M_\mu\inv(z),$ for any law $\sigma$ on $[0,+\infty)$, \be\label{1.7.8.3}M_{\sigma\boxtimes \mu_\la} \inv=\f{z}{z+1}\f{S_{\sigma}}{1+\la z}=\f{M_\sigma^{\lan -1\ran}}{1+\la z}\;\quad \textrm{ and }\;\quad M_{\sigma\bxs \mu_\la}\inv=(1+\la z)M_\sigma\inv.\ee At last, since $\bxp\delta_1=*\delta_1$, which implies that $M_{\sigma\bxp\delta_1}(z)=[(z+1)M_\sigma(z)+z]\circ \f{z}{z+1}$, for any symmetric law $\nu$, we have \be\label{1.07.08.1}M_{((\nu^2\ltimes\mu_\la)\bxp\delta_1)\boxtimes\mu_\la}\inv=\ff{1+\la z}\times\f{z}{1+z}\circ\lf[(z+1)\lf((1+\la z)M_{\nu^2}\inv\ri)\inv+z\ri]\inv .\ee

\subsection{Proof of the result}\label{2.7.8.3}
So let us consider a symmetric \pro measure $\nu$ \st $\nu^2\bxs \mu_\la$ exists and let us prove \eqref{30.06.08.2}. As proved in the proof of Theorem 3.8 of \cite{bg07}, for any symmetric \pro measure $\tau$, $H_\tau$ characterizes $\tau$, thus it suffices to prove that $H_{ \nu\arc \sqrt{\mu_\la}}=H_{m}$ for $m=\sqrt{ [(\nu^2\bxs\mu_\la)\bxp\delta_1]\bxt\mu_\la}.$ By Theorem 4.3 and the paragraph preceding in \cite{fbg05.inf.div},
 $C_{\sqrt{\mu_\la}}(z)=z$. 
Thus  Lemma 4.1 of \cite{bba07} applies here, and it states that  in a neighborhood of zero in $\C\bck [0,+\infty)$,    $$H_{\nu\arc \sqrt{\mu_\la}}=H_\nu\circ \lf(\f{H_\nu}{T(H_\nu+M_{\nu^2})}\ri)^{\lan-1\ran},$$where $T(z)=(\la z+1)(z+1)$.  So it suffices to prove that in such a neighborhood of zero, $$H_m=H_\nu\circ \lf(\f{H_\nu}{T(H_\nu+M_{\nu^2})}\ri)\inv,\quad\textrm{ i.e. }\quad H_m\circ \f{H_\nu}{T(H_\nu+M_{\nu^2})}=H_\nu.$$Using the fact that for any symmetric law $\tau$, $H_\tau(z)=zT(M_{\tau^2}(z)),$ it amounts to prove that $$\f{H_\nu}{T(H_\nu+M_{\nu^2})}\times T\circ M_{m^2}\circ \f{H_\nu}{T(H_\nu+M_{\nu^2})}  =H_{\nu}(z), $$i.e.$$
T\circ M_{m^2}\circ \f{H_\nu}{T(H_\nu+M_{\nu^2})}=T(H_\nu(z)+M_{\nu^2}(z)),$$which is implied, simplifying by $T$ and using again $H_\tau(z)=zT(M_{\tau^2}(z))$, by $$  M_{m^2}\circ \f{zT(M_{\nu^2}(z))}{T[zT(M_{\nu^2}(z))+M_{\nu^2}(z)]}=zT[M_{\nu^2}(z)]+M_{\nu^2}(z).$$ It is implied,  composing by  $M_{\nu^2}\inv$ on the right and by $M_{m^2}\inv$ on the left, by $$ M_{\nu^2}\inv\times T=(T\times M_{m^2}\inv)\circ (M_{\nu^2}\inv(z)T(z)+z).$$Using the expression of $M_{m^2}\inv$ given by \eqref{1.07.08.1}, it amounts to prove that $$ M_{\nu^2}\inv(z)T(z)=$$ $$
\lf((z+1)\times\f{z}{1+z}\circ\lf[(z+1)\lf((1+\la z)M_{\nu^2}^{\lan -1\ran}\ri)^{\lan -1\ran}+z\ri]\inv\ri)\circ (M_{\nu^2}\inv(z)T(z)+z),$$i.e. that
$$ \f{M_{\nu^2}\inv(z)T(z)}{M_{\nu^2}\inv(z)T(z)+z+1}=\f{z}{1+z}\circ\lf[(z+1)\lf((1+\la z)M_{\nu^2}^{\lan -1\ran} \ri)^{\lan -1\ran}+z\ri]\inv\circ (M_{\nu^2}\inv(z)T(z)+z)
.$$Now, composing by $\lf[(z+1)\lf((1+\la z)M_{\nu^2}^{\lan -1\ran} \ri)^{\lan -1\ran}+z\ri]\circ\f{z}{1-z}$ on the left, it gives
$$ \lf[(z+1)\lf((1+\la z)M_{\nu^2}^{\lan -1\ran} \ri)^{\lan -1\ran}+z\ri]\circ[(1+\la z)M_{\nu^2}\inv(z)]=M_{\nu^2}\inv(z)T(z)+z,
$$i.e.
$$ [M_{\nu^2}\inv(z)(\la z+1)+1]{z}+[M_{\nu^2}\inv(z)(\la z+1)]=M_{\nu^2}\inv(z)(\la z+1)(z+1)+z,
$$which is easily verified.

\subsection{Remarks on this result}
\begin{rmq}\label{1.7.8.1}{\rm Note that we did not use the fact that  $\nu^2\bxs \mu_\la$ exists  to prove that $H_{ \nu\arc \sqrt{\mu_\la}}=H_{m}$. It means that if $\nu^2\bxs \mu_\la$ doesn't exist,  there is no more \pro measure $\mu$ on $[0,+\infty)$ \st $M_\mu\inv=(1+\la z)M_{\nu^2}\inv$  as in \eqref{1.7.8.3},  but the polynomial expression of the  moments of $ \nu\arc \sqrt{\mu_\la}$ (i.e. of the limit symmetrized singular law of the matrix $A_n+G_n$ of Theorem \ref{1.07.08.3}) in the moments of $\nu$ following from $H_{ \nu\arc \sqrt{\mu_\la}}=H_{m}$ for $m=\sqrt{ [(\nu^2\bxs\mu_\la)\bxp\delta_1]\bxt\mu_\la}$  stays true (see Remark \ref{1.7.8.16h}).}\end{rmq}

\begin{rmq}[Case $\la=0$]\label{1.7.8.2}{\rm  A continuous way to define $\mu_\la$ for any $\la\in [0,1]$ is to define it to be the \pro measure with free cumulants $k_n(\mu_\la)=\la^{1-n}$ for all $n\geq 1$ (see \cite{ns06}). This definition gives $\mu_0=\delta_1$. Note that by definition of the rectangular free convolution with null ratio $\bxp_0$ (which is recalled in the introduction), the relation \eqref{30.06.08.2} stays true for $\la=0$.}\end{rmq}

\begin{rmq}\label{1.7.8.17h39}{\rm Note that the original proof of Debbah and Ryan in \cite{dr07} is based on the combinatorics approach to freeness, via the free cumulants of Nica and Speicher \cite{ns06}, whereas our proof is based on the analytical machinery for the computation of the rectangular $R$-transform, namely the rectangular $R$-transform. It happens sometimes that combinatorial proofs can be translated on the analytical plan by considering the generating functions of the combinatorial objects in question. Notice however that it is not what we did here. Indeed, the  rectangular $R$-transform machinery is  actually related to other cumulants than the ones of Nica and Speicher. These are the so-called rectangular cumulants, defined in \cite{bg07}.}\end{rmq}

\subsection{Remarks about the free deconvolution by $\mu_\la$}
The following corollary is part of the answer given in the present paper to the question asked in the last section of the paper of Debbah and Ryan \cite{dr07}. Let us endow the set of \pro measures on the real line with the weak topology \cite{billingsley}. 

\begin{cor}The functional $\nu\mapsto [(\nu^2\bxs\mu_\la)\bxp\delta_1]\bxt\mu_\la$, defined on the set of \pro measures $\nu$ on $[0,+\infty)$ \st $\nu\bxs\mu_\la$ exists, extends continuously to the whole set of \pro measures  on $[0,+\infty)$.\end{cor} 

\begin{pr} We just proved, in section \ref{2.7.8.3}, that the formula $$\nu\arc \sqrt{\mu_\la}=\sqrt{ [(\nu^2\bxs\mu_\la)\bxp\delta_1]\bxt\mu_\la}$$is true for any \pro measure $\nu$ on $[0,+\infty)$. 
Since the operation $\arc$ is  continuous on the set of symmetric \pro measures on the real line (Theorem 3.12 of \cite{bg07}) and the bijective corespondance between symmetric laws on the real line and laws on $[0,+\infty)$, which maps any symmetric law to its push-forward by the map $t\mapsto t^2$, is continuous with continuous inverse, the corollary is obvious.\end{pr}

 The  functional  $\bxs\mu_\la$, which domain is contained in the set  of \pro measures on $[0,+\infty)$, plays surprisingly a key role here. It seems natural to try to study its domain. The first step is to notice that this domain is the whole set of \pro measures on $[0,+\infty)$ if and only if $\delta_1$ is in this domain, and that in this case,  the functional $\bxs\mu_\la$ is simply equal to $\bxt(\delta_1\bxs\mu_\la)$.
However, the following proposition states that despite  the previous corollary, the domain of  the  functional  $\bxs\mu_\la$ is not the whole set of \pro measures on $[0,+\infty)$.
 
\begin{propo}The Dirac mass $\delta_1$ at $1$ is not in the domain of the functional  $\bxs\mu_\la$.\end{propo}
 
 \begin{pr} Suppose that there is a \pro measure $\tau$ on $[0,+\infty)$ \st $\delta_1=\tau\bxt\mu_\la$. 
Such a law $\tau$ has to satisfy $S_\tau(z)=1+\la z$. It implies that for $z$ small enough, $M_\tau(z)=\f{z-1+[(1-z)^2+4\la z]^{1/2}}{2\la}$. Such a function doesn't admit any analytic continuation to $\C\bck[0,+\infty)$, thus no such \pro measure $\tau$ exists.
  \end{pr}

\section{Relations between square and rectangular matrices/convolutions} 

The theorem of Debbah and Ryan gives an expression of the empirical singular 
measure of the sum of two rectangular random matrices in terms of operations related to hermitian square random matrices. Two other results   relate empirical singular measures of (non hermitian) square or rectangular random matrices to the operations devoted to hermitian random matrices.

The first one can be resumed by $\bxp_1=\bxp$. Concretely, it states that denoting by $\operatorname{ESM}(X)$ the symmetrization of empirical singular 
measure of any rectangular matrix $X$, for any pair $M,N$ of large $n$ by $p$ random matrices, one of them being invariant in law by the left and right actions of the unitary groups, for $n/p\simeq 1$, \be\label{2.7.8.16h008}\esm(M+N)\simeq \esm(M)\bxp \esm(N). \ee Note that the matrices $M,N$ are not hermitian, which makes  \eqref{2.7.8.16h008} pretty surprising (since $\bxp$ was defined with hermitian random matrices). It means that for  $\eps, \eps_1,\eps_2$ independent random variables with law $\f{\delta_{-1}+\delta_1}{2}$                     independent of $M$ and $N$, we have \be\label{1.7.8.20h14}\operatorname{Spectrum}(\eps |M+N|){\simeq}\operatorname{Spectrum}(\eps_1 |M|+\eps_2|N|)\ee

The second  one can be resumed by for any pair $\nu,\tau$ of symmetric \pro measures on the real line, $(\nu\bxp_0\tau)^2=\nu^2\bxp\tau^2$
 Concretely, it states that for any pair $M,N$ of $n$ by $p$ unitarily invariant random matrices,  for $1<<n<<p$, 
\be\label{1.7.8.20h15}\operatorname{Spectrum}[(M+N)(M+N)^*] {\simeq}\operatorname{Spectrum}( MM^*+NN^*).\ee

The advantage of the result of Debbah and Ryan on those ones is that it works for any value of the ratio $\la$, but its disadvantage is that it only works when one of the laws convoluted is $\mu_\la$, i.e.  one of the matrices considered is a Gaussian one. In fact this sharp restriction can be understood by the fact that among rectangular random matrices which are invariant in law under multiplication by unitary matrices, the Gaussian ones are the only ones which can be extended to square matrices which are also invariant in law under multiplication by unitary matrices. 

It could be interesting to understand better how relations like \eqref{1.7.8.20h14}, \eqref{1.7.8.20h15} or  like the one of the Debbah and Ryan's theorem work and can be generalized. Unfortunately, until now, even though   nice proofs (see \cite{bg07} for \eqref{1.7.8.20h14} and \eqref{1.7.8.20h15} or Theorem 4.3.11 of \cite{hiai} and Proposition 3.5 of \cite{haag2} for the $n=p$ case of \eqref{1.7.8.20h14}) relying in free \pro have been given for these results relating rectangular convolutions  and ``square non hermitian convolutions" with the ``square hermitian convolution" (i.e. $\bxp$), no ``concrete" explanation has been given, and no generalization (to any $\la$, to any pair of \pro measures) neither.
Such a  generalization could be  the given of a functional $f_\la$ on the set of symmetric \pro measures \st for all $\nu,\tau$ symmetric \pro measures, $\nu\arc\tau$ is the only symmetric \pro measure satisfying $$f_\la(\nu\arc\tau)=f_\la(\nu)\bxp f_\la(\tau).$$ 
 Note that in the case $\la=1$, the functional $f_\la(\nu)=\nu$ works, and in the case $\la=0$, the functional which maps a measure to its push-forward by the square function works.

\begin{rmq}{\rm  Let $(\A,\vfi)$ be  a $*$-non commutative \pro space and $p_1, p_2$ be two self-adjoint projectors of $\A$ \st $p_1+p_2=1$ \st $\la=\vfi(p_1)/\vfi(p_2)$. %Then there is a natural $*$-algebra identification  which allows us to identify the elements of $\A$ with two by two matrices which  $i,j$-th entry belongs to $p_i\A p_j$, for all $i,j\in\{1,2\}$.  
 As explained in Proposition-Definition 2.1 of \cite{bg07}, $\arc$ can be defined by the fact that for any pair $a,b\in p_1\A p_2$ free with amalgamation over  $\operatorname{Vect}(p_1,p_2)$, the symmetrized   distribution of $|a+b|$ in $(p_1\A p_1,\ff{\vfi(p_1)}\vfi_{|p_1\A p_1})$ is the rectangular free convolution with ratio $\la$ of the symmetrized   distributions of $|a|$ and $|b|$ in the same space.

Moreover, it is easy to see that   for all $a\in p_1\A p_2$, the symmetrized distribution $\tau$  of   $|a|$ in $(p_1\A p_1,\ff{\vfi(p_1)}\vfi_{|p_1\A p_1})$ is linked to  the distribution $\nu$ of  $a+a^*$ %$\begin{bmatrix}0&a\\ a^*&0\end{bmatrix}$
 in $(\A,\vfi)$ by the relation  $\nu= \f{2\la}{1+\la}\tau+\f{1-\la}{1+\la}\delta_0.$%where $\la=\vfi(p_1)/\vfi(p_2)$. 

When $\la=1$, the equation $\bxp=\arc$ can be summurized in the following way: for $a,b\in p_1\A p_2$  free with amalgamation over $\operatorname{Vect}(p_1,p_2)$, the distribution of $(a+b)+(a+b)^*$ in $(\A, \vfi)$ is the free convolution of the distributions of  $a+a^*$ and $b+b^*$.
%$$\begin{bmatrix}0&a_1+a_2\\ a_1^*+a_2^*&0\end{bmatrix} $$ is the free convolution of the distributions of  $\begin{bmatrix}0&a_1\\ a_1^*&0\end{bmatrix}$ and $\begin{bmatrix}0&a_2\\ a_2^*&0\end{bmatrix}$.

If this had stayed true for other values of $\la$, 
it would have meant that for all $\nu,\tau$ compactly supported symmetric probability measures on the real line, we have \begin{equation}\label{13.03.06.1}f_\la(\nu\arc\tau)=f_\la(\nu)\bxp f_\la(\tau) ,\end{equation} where $f_\la$ is the function which maps a \pro measure $\tau$ on the real line  to $\f{2\la}{1+\la}\tau+\f{1-\la}{1+\la}\delta_0$. But looking at fourth moment, it appears that \eqref{13.03.06.1} isn't true.}\end{rmq}


\begin{thebibliography}{99}\bibitem[BBA07]{bba07} Belinschi, S., Benaych-Georges, F., Guionnet, A. \emph{Regularization by free additive convolution, square and rectangular cases}. To appear in {\bf Complex Analysis and Operator Theory}. 
\bibitem[BG07a]{fbg05.inf.div} Benaych-Georges, F. \emph{Infinitely divisible distributions for rectangular free convolution: classification and matricial interpretation}
{\bf Probability Theory and Related Fields} Volume 139, Numbers 1-2 / septembre 2007, 143-189. 
\bibitem[BG07b]{bg07} Benaych-Georges, F. \emph{Rectangular random matrices, related convolution}. To appear in {\bf Probability and Theory Related Fields}.
\bibitem[B68]{billingsley} Billingsley, P. \emph{Convergence of probability measures} Wiley, 1968
\bibitem[DR07]{dr07} Debbah, M., Ryan, {\O}. {\em Multiplicative free Convolution and Information-Plus-Noise Type Matrices}. arXiv.
\bibitem[DR08]{dr08} Debbah, M., Ryan, {\O}. {\em Free Deconvolution for Signal Processing Applications} Second round review, submitted to  IEEE
transactions on Information Theory, 2007.
\bibitem[HL00]{haag2}  Haagerup, U.,  Larsen, F. \emph{Brown's spectral distribution measure for R-diagonal elements in finite von Neumann algebras} {\bf Journ. Functional Analysis} 176, 331-367 (2000).
\bibitem[HP00]{hiai} Hiai, F., Petz, D. \emph{The semicircle law, free random variables, and entropy} Amer.
Math. Soc., Mathematical Surveys and Monographs Volume 77, 2000
\bibitem[NS06]{ns06} Nica, Alexandru; Speicher, Roland {\em Lectures on the combinatorics of free probability}. London Mathematical Society Lecture Note Series, 335. Cambridge University Press, Cambridge, 2006.
\bibitem[VDN91]{vdn91} Voiculescu, D.V., Dykema, K., Nica, A. \emph{Free random variables} CRM Monograghs Series No.1, Amer. Math. Soc., Providence, RI, 1992 
\end{thebibliography}
\end{document}